\documentclass[12pt]{article}

\newcommand{\End}{\mathop{\mathrm{End}}}
\newcommand{\Gr}{\mathop{\mathrm{Gr}}}
\newcommand{\Hull}{\mathop{\mathrm{Hull}}}

\newtheorem{theorem}{Theorem}[section]
\newtheorem{cor}[theorem]{Corollary}
\newtheorem{prop}[theorem]{Proposition}

\newtheorem{unnumber}{}

\newtheorem{prepf}[unnumber]{Proof}
\newenvironment{pf}{\prepf\rm}{\endprepf}

\begin{document}

\title{Dixon's Theorem and random synchronization}
\author{Peter J. Cameron\\School of Mathematical Sciences\\
Queen Mary, University of London\\Mile End Road\\London E1 4NS, U.K.}
\date{}
\maketitle

\begin{abstract}
A transformation monoid on a set $\Omega$ is called \emph{synchronizing} if
it contains an element of rank~$1$ (that is, mapping the whole of $\Omega$
to a single point). In this paper, I tackle the question: given $n$ and $k$,
what is the probability that the submonoid of the full transformation monoid
$T_n$ generated by $k$ random transformations is synchronizing?

The question has some similarities with a similar question about the probability
that the subgroup of $S_n$ generated by $k$ random permutations is transitive.
For $k=1$, the answer is $1/n$; for $k=2$, Dixon's Theorem asserts that it is
$1-o(1)$ as $n\to\infty$ (and good estimates are now known). For our
synchronization question, for $k=1$ the answer is also $1/n$; I conjecture
that for $k=2$ it is also $1-o(1)$.

Following the technique of Dixon's theorem, we need to analyse the maximal
non-synchronizing submonoids of $T_n$. I develop a very close connection 
between transformation monoids and graphs, from which we obtain a description
of non-synchronizing monoids as endomorphism monoids of graphs satisfying some
very strong conditions. However, counting such graphs,
and dealing with the intersections of their endomorphism monoids, seems
difficult.

\paragraph{Keywords:} transformation monoid, synchronization, graph
homomorphisms, random generation.
\end{abstract}

\section{Dixon's Theorem}

In 1969, John Dixon~\cite{dixon} proved the following theorem:

\begin{theorem}
The probability that two random permutations in the symmetric group $S_n$
generate $S_n$ or $A_n$ is $1-o(1)$ as $n\to\infty$.
\end{theorem}

In fact, good estimates are known. Babai~\cite{babai} showed that the
probability is $1-1/n+O(1/n^2)$: the term $1/n$ arises from the probability
that the two permutations have a common fixed point. Several further terms
of the asymptotic expansion are known.

It is my purpose here to begin a similar analysis for the \emph{full
transformation monoid} $T_n$ on the set $\Omega=\{1,\ldots,n\}$. Things
are a little different, since $T_n$ requires three generators. (If the
monoid $M$ is generated by a set $S$ of transformations, then the group
of permutations in $M$ is generated by the permutations in $S$; so
if $M=T_n$ with $n>2$, then $S$ must contain at least two permutations,
and at least one non-permutation.) Indeed, since permutations are exponentially
scarce in $T_n$, we have to choose a huge number of random elements in order
to generate $T_n$ with high probability.

Further analysis of Dixon's theorem suggests a different approach. The
first, and easier, step is to calculate the probability that two
permutations in $S_n$ generate a transitive subgroup. If $c_n$ is the
number of pairs of elements of $S_n$ which generate a transitive subgroup,
then counting pairs according to the orbit of the point $1$ of the group
they generate gives
\[\sum_{k=1}^n{n-1\choose k-1}c_k((n-k)!)^2 = (n!)^2,\]
a recurrence relation from which $c_n$ can be determined. It is then easy
to show that $c_n/(n!)^2 = 1-1/n+O(1/n^2)$.

However, a cruder analysis is more useful in other situations. The
maximal intransitive subgroups of $S_n$ have the form $S_k\times S_{n-k}$
for $1\le k\le\lfloor n/2\rfloor$. If two elements fail to generate a
transitive subgroup, then they lie in some maximal intransitive subgroup;
the probability of this is at most
\[\frac{1}{(n!)^2}\sum_{k=1}^{\lfloor n/2\rfloor}{n\choose k}(k!)^2
((n-k)!)^2 = \frac{1}{n}+O\left(\frac{1}{n^2}\right).\]

The remainder of the proof of Dixon's Theorem involves showing that the
probability that the two permutations lie in a transitive subgroup other
than the symmetric or alterating group is very small. This probability is
estimated similarly to the above, by bounding the number and order of maximal
transitive subgroups other than $S_n$ and $A_n$.

We note in passing that the probability that a single random permutation
in $S_n$ generates a transitive subgroup is $1/n$. For the permtations which
generate transitive subgroups are the $n$-cycles, and it is well-known that
there are exactly $(n-1)!$ of these.

\section{Synchronizing monoids}

Let $T_n$ be the full transformation monoid on the set $\Omega=\{1,\ldots,n\}$,
consisting of all \emph{endofunctions} $f:\Omega\to\Omega$, with the operation
of composition. A \emph{transformation monoid} is a submonoid of $T_n$.

A transformation monoid $M$ is said to be \emph{synchronizing} if it contains
an element of rank~$1$ (that is, a function whose image has cardinality~$1$).
It seems that synchronizing monoids behave a little like transitive subgroups
of $S_n$. The first observation gives an exact parallel:

\begin{prop}
The probability that a random endofunction generates a synchronizing monoid
is $1/n$.
\end{prop}

\begin{pf}
The endofunction $f$ generates a synchronizing monoid if and only if it has
a unique periodic point. Such a function is defined by a rooted tree, with
edges directed towards the root. There are $n^{n-1}$ rooted trees, and $n^n$
endofunctions altogether.
\end{pf}

I conjecture that the probability that two random endofunctions generate a
synchronizing monoid is $1-o(1)$. The strategy is to describe the maximal
non-synchronizing monoids, and then to argue as in the proof of Dixon's
theorem. The first part of the programme is realised here, and some evidence
towards the second is given.

Of course, the analogy between transitive subgroups and synchronizing 
submonoids is not perfect. {\L}uczak and Pyber~\cite{lp} showed that the
proportion of elements of $S_n$ which lie in transitive subgroup of $S_n$
except $S_n$ and possibly $A_n$ is $1-o(1)$, though the rate of convergence
is not well understood. However, every element of $T_n$ lies in a proper
synchronizing submonoid. For, if $\langle g\rangle$ is synchronizing, then
so is $\langle f,g\rangle$ for any $f\in T_n$; but $\langle f,g\rangle
\neq T_n$, since $T_n$ requires at least three generators.

\section{Monoids and graphs}

There is a very close connection between transformation monoids and graphs
which we define in this section. It has some features of a Galois
correspondence, but things are not quite so simple.

Let $\Omega=\{1,\ldots,n\}$. We define maps in each direction between 
transformation monoids on $\Omega$ and graphs on the vertex set $\Omega$.

One direction is well-known. Given a graph $X$, an \emph{endomorphism} of $X$
is an endofunction on $\Omega$ which maps edges of $X$ to edges. (We do not
care what it does to non-edges, which may be mapped to non-edges or to edges
or to single vertices). The endomorphisms of $X$ clearly form a monoid
$\End(X)$.

In the other direction, given a transformation monoid $M$, we define a graph
$X=\Gr(M)$ by the rule that two vertices $v,w$ are adjacent if and only if
there does not exist $f\in M$ such that $vf=wf$.

Not every graph occurs as the graph of a monoid. Recall that the \emph{clique
number} $\omega(X)$ is the cardinality of the largest complete subgraph of
$X$, and the \emph{chromatic number} $\chi(X)$ is the smallest number of
colors required for a proper colouring of the vertices (so that adjacent
vertices have different colours). Clearly $\omega(X)\le\chi(X)$, since
all vertices in a clique must have different colours; these parameters may
differ arbitrarily.

\begin{theorem}
For any transformation monoid $M$, $\omega(\Gr(M))=\chi(\Gr(M))$, and this
number is equal to the minimum rank of an element of $M$.
\label{t:cc}
\end{theorem}

\begin{pf}
Let $f$ be an element of $M$ of minimum rank, and let $S$ be the image of $f$.
Then the induced subgraph on $S$ is a clique; for if $v,w\in S$ are not
adjacent, then there exists $g\in M$ with $vg=wg$, so that $fg$ has smaller
rank than $f$. But the map $f$ is a proper colouring of $\Gr(M)$, since by
definition the images of adjacent vertices are distinct. So we have
$\chi(\Gr(M))\le|S|\le\omega(\Gr(M))$, whence equality holds throughout.
\end{pf}

\begin{cor}
\begin{enumerate}
\item $\Gr(M)$ is a complete graph if and only if $M\le S_n$ (that is, all
elements of $M$ are permutations).
\item $\Gr(M)$ is a null graph if and only if $M$ is synchronizing.
\item If $M_1\le M_2$, then $\Gr(M_2)$ is a spanning subgraph of $\Gr(M_1)$.
\end{enumerate}
\label{c:cn}
\end{cor}

The first and third parts, and the reverse implication in the second, are
clear; the forward implication in the second part follows immediately from the
preceding Theorem.

The two maps (from graphs to monoids and from monoids to graphs) are not
mutually inverse, and do not (quite) form a Galois connection; but they do
satisfy the following:

\begin{theorem}
For any transformation monoid $M$,
\begin{enumerate}
\item $M\le\End(\Gr(M))$;
\item $\Gr(\End(\Gr(M)))=\Gr(M)$.
\end{enumerate}
\label{t:gm}
\end{theorem}

\begin{pf}
(a) Let $f\in M$, and let $\{v,w\}$ be an edge in $\Gr(M)$. By definition,
$vf\ne wf$. Could $vf$ and $wf$ be non-adjacent in $\Gr(M)$? If so, then there
would be $g\in M$ such that $(vf)g=(wf)g$. But then the map $fg\in M$
satisfies $v(fg)=w(fg)$, contradicting the fact that $v$ and $w$ are joined.
So $f\in\End(\Gr(M))$.

(b) If $\{v,w\}$ is an edge of $\Gr(M)$, then no endomorphism of $\Gr(M)$
collapses it to a point, and so $\{v,w\}$ is an edge of $\Gr(\End(\Gr(M))$.
Conversely, suppose that $v$ and $w$ are not adjacent in $\Gr(M)$. Then by
definition there exists $f\in M$ such that $vf=wf$. Since $f\in\End(\Gr(M))$
by (a), we see tat $v$ and $w$ are not adjacent in $\Gr(\End(\Gr(M)))$. So
these two graphs are equal.
\end{pf}

Given a graph $X$, the graph $\Gr(\End(X))$ is called the \emph{hull} of $X$,
and is studied in \cite{ck}. Theorem~\ref{t:gm}(b) shows that 
$\Hull(\Hull(X))=\Hull(X)$. In other words, a graph $X$ is a hull if and only
if it is its own hull (that is, $\Hull(X)=X$).

\section{Another construction}

Here is another construction which doesn't decrease the endomorphism monoid
of a graph.

\begin{prop}
Let $X$ be a graph on the vertex set $\Omega$ with $\omega(X)=m$.
Let $X'$ be the spanning subgraph of $X$ which consists of those edges
of $X$ which are contained in cliques of size $m$. Then
$\End(X)\le\End(X')$.
\label{p:deriv}
\end{prop}

\begin{pf}
Suppose not. Then there exists $f\in\End(X)$ such that $f\notin\End(X')$. This
means that there is an edge $\{v,w\}$ of $X'$ such that either $vf=wf$,
or $\{vf,wf\}$ is a non-edge of $X'$.

The first case is impossible since $\{v,w\}$ is an edge of $X$ and
$f\in\End(X)$. Suppose that the second case happens. Then $\{vf,wf\}$ is an
edge of $X$, and was deleted because it is not contained in any clique of
size $m$. But $\{v,w\}$ is not deleted, so lies in a clique $C$ of
$X$ with $|C|=m$; and then $Cf$ is a clique of $X$ with
$\{vf,wf\}\subseteq Cf$ and $|Cf|=m$, a contradiction.
\end{pf}

I will call $Y$ the \emph{derived graph} of $X$.

\section{Maximal non-synchronizing monoids}

In this section we will give a description of the maximal non-synchronizing
monoids in terms of graphs. Note that, if the graph $X$ is non-null, then
$\End(X)$ is non-synchronizing. The main theorem is the following:

\begin{theorem}
Let $M$ be a maximal non-synchronizing submonoid of $T_n$. Then there are
graphs $X$ and $Y$ on the vertex set $\Omega=\{1,\ldots,n\}$ satisfying the
following conditions:
\begin{enumerate}
\item $\End(X)=\End(Y)=M$;
\item $\omega(X)=\omega(Y)=\chi(X)=\chi(Y)$;
\item $X=\Hull(Y)$;
\item $Y=X'$.
\end{enumerate}
\label{t:mns}
\end{theorem}

\begin{pf}
Let $M$ be maximal non-synchronizing. Let $X=\Gr(M)$ and $Y=X'$.
Then $X$ has at least one edge (by Corollary~\ref{c:cn}(b)), and satisfies
$\omega(X)=\chi(X)$ (by Theorem~\ref{t:cc}). Moreover, $M\le\End(X)$, by
Theorem~\ref{t:gm}(a); maximality of $M$ implies that equality holds.

Now $M=\End(X)\le\End(Y)$ by Proposition~\ref{p:deriv}; maximality of $M$
implies that equality holds. Furthermore, it is clear that
\[\omega(Y)=\omega(X)=\chi(X)\ge\chi(Y)\ge\omega(Y),\]
so equality holds throughout. Finally, since $\End(X)=\End(Y)$, we see that
\[X=\Hull(X)=\Gr(\End(X))=\Gr(\End(Y))=\Hull(Y).\]
\end{pf}

I do not know any examples where $X$ and $Y$ are not equal. If they are
equal, then the converse holds:

\begin{theorem}
Let $X$ be a hull (other than the null graph), in which every edge is
contained in a clique of size $\omega(X)$. Then $\End(X)$ is a maximal
non-synchronizing submonoid of $T(\Omega)$.
\label{t:nearcon}
\end{theorem}

\begin{pf}
Let $f$ be any endofunction not in $M=\End(X)$. By Corollary~\ref{c:cn}(b),
it suffices to show that for any $v,w\in\Omega$, there is an element
$g\in M'=\langle M,f\rangle$ such that $vg=wg$. Since $X$ is a hull, this holds
for any $v,w$ for which $\{v,w\}$ is a non-edge of $X$, so we may assume
that $\{v,w\}$ is an edge.

I claim that, if $\{v',w'\}$ is another edge, then there is an endomorphism
$h$ of $X$ satisfying $vh=v'$ and $wh=w'$. For, by assumption, there is a
clique $C$ with $|C|=\omega(X)$ containing $v'$ and $w'$; now there is an
endomorphism from $X$ onto $C$, and since $C$ is complete, we may order its
elements arbitrarily, so that in particular the images of $v$ and $w$ are
$v'$ and $w'$ as claimed.

Since $f$ is not an endomorphism, there is an edge $\{x,y\}$ of $X$ such
that either $xf=yf$, or $\{xf,yf\}$ is a non-edge. Composing $f$ with an
endomorphism if necessary, we may assume that $xf=yf$. Taking $v'=x$ and
$w'=y$, and composing $h$ of the preceding paragraph with $f$, we find an
element of $M$ with the reqired property.
\end{pf}

There are many graphs satisfying the hypotheses of this theorem. The smallest
consists of a single edge; there are $n(n-1)/2$ graphs of this form and
each has $2n^{n-2}$ endomorphisms. So the probability that a random pair of
endofunctions are both endomorphisms of a graph of this form is at most
\[\frac{n(n-1)}{2}\,\frac{n^{2(n-2)}}{n^{2n}} = O(n^{-2}).\]

This suggests that the probability that two random endofunctions generate a
synchronizing monoid is at least $1-O(1/n^2)$. However, we are still some
way from a proof, since there are many graphs that need to be considered.
Of course, there are big overlaps between their endomorphism monoids, so
inclusion-exclusion will have to be applied much more carefully than in the
case of Dixon's Theorem.

\section{Open problems}

The main problem is to prove that the probability that two random elements
generate a synchronizing monoid is $1-o(1)$.

A variant is to choose $r+s$ elements, of which $r$ are random permutations
and the remaining $s$ are random endofunctions. If $r\ge2$ and $s\ge1$, then
by Dixon's Theorem the permutations generate $S_n$ or $A_n$ with high
probability, and the entire monoid is synchronizing with
high probability. The interesting case here is $r=s=1$.

A final problem is whether the two graphs in Theorem~\ref{t:mns} can be
distinct. If not, then the conditions of Theorem~\ref{t:nearcon} would be
necessary and sufficient for a monoid to be maximal non-synchronizing.

\end{document}